\documentclass[11pt]{amsart}
\usepackage{amssymb}

\newtheorem{lemma}{Lemma}[section]

\newtheorem{theorem}[lemma]{Theorem}
\newtheorem{remark}[lemma]{Remark}
\newtheorem{proposition}[lemma]{Proposition}
\newtheorem{corollary}[lemma]{Corollary}
\newtheorem{definition}[lemma]{Definition}
\newtheorem{example}[lemma]{Example}
\newtheorem{conjecture}[lemma]{Conjecture}

\begin{document}
\title{Topology of real Milnor fibrations for non-isolated singularities}
\author{Nicolas Dutertre and Raimundo Ara\'ujo dos Santos}

\address{N. Dutertre: Aix-Marseille Universit\'e, LATP, 39 rue Joliot-Curie, 13453 Marseille Cedex 13, France}

\email{dutertre@cmi.univ-mrs.fr}

\address{R. N. Ara\'ujo dos Santos: Departamento de Matem\'{a}tica, Instituto de Ci\^{e}ncias Matem\'{a}ticas e de Computa\c{c}\~{a}o, Universidade de S\~{a}o Paulo - Campus de S\~{a}o Carlos, Caixa Postal 668, 13560-970 
S\~ao Carlos, SP, Brazil}

\email{rnonato@icmc.usp.br}


\begin{abstract}
We consider a real analytic map $F=(f_1,\ldots,f_k) : (\mathbb{R}^n,0) \rightarrow (\mathbb{R}^k,0)$, $2 \le k \le n-1$,  that satisfies Milnor's conditions (a) and (b) introduced by D. Massey.
This implies that every real analytic $f_I=(f_{i_1},\ldots,f_{i_l}) : (\mathbb{R}^n,0) \rightarrow (\mathbb{R}^l,0)$, induced from $F$ by projections where $1 \le l \le n-2$ and $I=\{i_1,\ldots,i_l\}$, also satisfies Milnor's conditions (a) and (b). We give several relations between the Euler characteristics of the Milnor fibre of $F$, the Milnor fibres of the maps $f_I$, the link of $F^{-1}(0)$ and the links of $f_I^{-1}(0)$.
\end{abstract}

\maketitle
\markboth{Nicolas Dutertre and Raimundo Ara\'ujo dos Santos}{Topology of real Milnor fibrations for non-isolated singularities}

\section{Introduction}

In the last fourty years many research has been developed toward understanding the geometry and topology of complex and real singularities. After the famous book of J. Milnor \cite{Mi} the search for real and complex topological invariants of singularity have gotten special attentions. In \cite{Mi} Milnor considered a holomorphic function $f:0\in U\subset \mathbb{C}^{n}\to \mathbb{C},$ $f(0)=0$ and $\nabla f(0)=0,$  and proved the existence of a smooth fiber bundle in a neighborhood of the critical point $0.$ Moreover, if the critical point is isolated he related the Euler-Poincar\'e number of the fiber of this fibration with the topological degree of the gradient vector field $\nabla f.$ This became a starting point of several others formulae in the real and complex settings. Let us remind below some of them in the real case.

\vspace{0.3cm}

In \cite{Kh} Khimshiashvili proved a Poincar\'e-Hopf formula which relates the Euler-Poincar\'e number of a regular local fiber of an  analytic function with the topological degree of its gradient vector field, as follows.

\vspace{0.2cm}

Let  $f : (\mathbb{R}^n,0) \rightarrow (\mathbb{R},0)$ be a germ of a real analytic function with isolated critical point, then
$$ \chi \big(f^{-1}(\delta) \cap B_\epsilon \big)=1-\hbox{sign}(-\delta)^n  \hbox{deg}_0 \nabla f,$$
where $0 < \vert \delta \vert \ll \varepsilon \ll 1$ is a regular value, $B_{\epsilon}$ stands for the close ball centered at the origin, $\nabla f$ is the gradient vector field of $f$ and deg$_0 \nabla f$ is the topological degree of the mapping
$$\epsilon \frac{\nabla f}{\| \nabla f \|} :S_\epsilon ^{n-1}  \rightarrow S_{\epsilon}^{n-1}.$$

A similar formula was proved before by Milnor (see \cite{Mi}, page 61) for the case of holomorphic functions with isolated singularity.

\vspace{0.3cm}

Several relative versions of the Khimshiashvili formula were proved afterwards. Let us present them briefly.
Let $\psi =(f_1,\ldots,f_k) : (\mathbb{R}^n,0) \rightarrow (\mathbb{R}^k,0)$, $2 \le k \le n$, be an analytic map-germ and let us denote by $\phi$ the map-germ $(f_1,\ldots,f_{k-1}): (\mathbb{R}^n,0) \rightarrow (\mathbb{R}^{k-1},0)$. We assume that $\psi^{-1}(0)$ and $\phi^{-1}(0)$ have an isolated singularity at $0$ (note that here the maps $\phi$ and $\psi$ do not need to have an isolated critical point at the origin). Some authors investigated the following difference:
$$D_{\delta,\alpha}= \chi \big( \phi^{-1}(\delta) \cap \{f_{k} \ge \alpha\} \cap B_\epsilon \big) - \chi \big( \phi^{-1}(\delta) \cap \{f_{k} \le \alpha\} \cap B_\epsilon \big),$$
where $(\delta,\alpha)$ is a regular value of $\psi$ such that $0 \le \vert \alpha \vert \ll \vert \delta \vert \ll \epsilon$.

\vspace{0.2cm}

In \cite{Du2}, the first author proved that:
$$D_{\delta,\alpha} \equiv \hbox{dim}_{\mathbb{R}} \frac{\mathcal{O}_{\mathbb{R}^n,0}}{I} \bmod 2,$$
where $\mathcal{O}_{\mathbb{R}^n,0}$ is the ring of analytic function-germs at the origin and $I$ is the ideal generated by $f_1,\ldots,f_{k-1}$ and all the $k\times k$ minors $\frac{\partial (f_k,f_1,\ldots,f_{k-1})}{\partial(x_{i_1},\ldots,x_{i_k})}$. This is only a mod 2 relation and we may ask if it is possible to get a more precise relation.

When $k=n$ and $f_k=x_1^2+\cdots+x_n^2$, according to Aoki et al. (\cite{AFN1}, \cite{AFS}), $D_{\delta,0}= \chi \big(\phi^{-1}(\delta) \cap B_\varepsilon \big)= 2 \hbox{deg}_0 H$ and $2\hbox{deg}_0 H$ is the number of semi-branches of $\phi^{-1}(0),$ where
$$H=(\frac{\partial (f_n,f_1,\ldots,f_{n-1})}{\partial(x_{1},\ldots,x_{n})},f_1,\ldots,f_{n-1}).$$
They proved a similar formula in the case $f_k=x_n$ in \cite{AFN2} and Szafraniec generalized all these results to any $f_k$ in \cite{Sz3}.

\vspace{0.2cm}

When $k=2$ and $f_2=x_1$, Fukui \cite{Fu1} stated that $$D_{\delta,0}=-\hbox{sign}(-\delta)^n \hbox{deg}_{0} H,$$ where $H=(f_1,\frac{\partial f_1}{\partial x_2},\ldots,\frac{\partial f_1}{\partial x_n})$. Several generalizations of Fukui's formula are given in \cite{Fu2}, \cite{Du1}, \cite{FK} and \cite{Du4}. Note that in \cite{Du4}, the first author gave degree formulas for the Euler characteristic of regular fibers of some map-germs from $(\mathbb{R}^n,0)$ to $(\mathbb{R}^2,0)$ called partially parallelizable.

\vspace{0.3cm}

More recently in \cite{ADD} the authors of the present paper and D. Dreibelbis proved an extension of Khimshiashvili's formula in the following way.

\vspace{0.3cm}

Take $\psi =(f_1,\ldots,f_k) : (\mathbb{R}^n,0) \rightarrow (\mathbb{R}^k,0)$, $n \geq k \geq 2$, an analytic map germ and suppose that $0\in \mathbb{R}^{n}$ is an isolated singular point of $\psi .$ By Milnor \cite{Mi}, page 98, it is known that for each small enough $\epsilon >0$, there exists $0<\eta \ll \epsilon $ such that the mapping

\begin{equation}\label{equation1}
\psi_{|}:B_{\epsilon}\cap \psi^{-1}(S_{\eta}^{p-1})\to S_{\eta}^{p-1}
\end{equation}
is the projection of a smooth locally trivial fibration.

\vspace{0.2cm}

Denote by $M_{\psi}$ its fiber, {\it also known as Milnor fiber}. The main result was:

\begin{proposition}[\cite{ADD}, page 71] Let $\psi =(f_1,\ldots,f_k) : (\mathbb{R}^n,0) \rightarrow (\mathbb{R}^k,0)$ as above.

\begin{itemize}

\item[(i)] If $n$ is even, then $\chi \big(M_{\psi}) = 1-\hbox{\em deg}_0 \nabla f_1.$ Moreover, $$\deg_{0} \nabla f_{1}=\cdots =\deg_{0} \nabla f_{k}.$$

\item[(ii)] If $n$ is odd, then $\chi \big(M_{\psi}) = 1 ,$ and $\hbox{\em deg}_0 \nabla f_i=0$, for all $i.$

\end{itemize}

\end{proposition}

\vspace{0.3cm}

In the present paper we will use tools from Morse theory and singularity theory to prove several extensions of the previous formula for the setting of real analytic mappings with non-isolated singularity. Our formulae relates the Euler-Poincar\'e number of the Milnor fibers with the Euler-Poincar\'e number of the singular link. We will assume that the analytic maps satisfies the Milnor conditions $(a)$ and $(b)$ defined by D. Massey in \cite{Ma}, so it implies the existence of {\it the tube Milnor fibration} like in ${\bf (1)}$ above. These conditions seem to be strong in the setting of real analytic mappings, but it is not difficult to show that any holomorphic function satisfies them. See Example \ref{R1} for further details. Therefore, our formulae also provide an extension of Milnor's formula (see \cite{Mi}, page 64) for the case of holomorphic functions with non-isolated singularity. We will also answer a question stated by Milnor in \cite{Mi}, page 100, (see below) under these more general Milnor conditions $(a)$ and $(b)$. Let us remind this conjecture below:

\vspace{0.3cm}

{\it `` Note that any polynomial mapping $\mathbb{R}^{n}\to \mathbb{R}^{k}$ with isolated singularity at origin can be composed with the projection $\mathbb{R}^{k}\to \mathbb{R}^{k-1}$ to obtain a new mapping $\mathbb{R}^{n}\to \mathbb{R}^{k-1},$ also with isolated singularity at origin "}.

\vspace{0.2cm}

\begin{conjecture}[\cite{Mi}, page 100]

The fiber of the fibration associated with this new mapping is homeomorphic to the product of the old fiber with the unit interval.
\end{conjecture}

We should say that as far as we know this problem was approached by A. Jacquemard in \cite{Ja} under different hypotheses that cover the isolated singular case as our hypotheses do.

\vspace{0.2cm}

The paper is organized as follows. In section 2, we remind the definition of Milnor's conditions $(a)$ and $(b)$ and recall the proof of Milnor's fibration theorem. Section 3 contains some auxiliary lemmas about subanalytic sets. In section 4, we state basic results for mappings satisfying Milnor's conditions. In section 5, we study the behavior of some critical points on the boundary of the Milnor fiber. In section 6, we consider a mapping satisfying Milnor's conditions $(a)$ and $(b),$  give a proof of Milnor's conjecture stated above and study the Euler-Poincar\'e numbers of the Milnor fibers of the mappings given by compositions of this initial mapping and projections. In section 7, we still consider these mappings and we relate the Euler-Poincar\'e number of their Milnor fibers to the Euler-Poincar\'e number of the links of their zero sets. In section 8, we establish several formulae for the Euler-Poincar\'e number of semi-analytic sets defined by the components of the initial mapping. Last section contains some applications and examples.

\vspace{0.2cm}

{\it Acknowledgments.} Computations are made in section 9, they have been carried out by a program, based on the Eisenbud-Levine-Khimshiashvili formula, written by A. \L ecki. The authors are very grateful to him and Z. Szafraniec for giving them this program. 

The authors thank the USP-Cofecub project ``UcMa133/12 - Structure fibr\'ee de l'espace au voisinage des singularités des applications". 

The first author is partially supported by the program ``Cat\'edras L\'evi-Strauss $-$ USP/French Embassy, no. 2012.1.62.55.7". 

\section{Milnor's conditions $(a)$ and $(b)$}

In this section we will follow the definitions and results given by D. Massey in \cite{Ma}.

Let $F=(f_1,\ldots,f_k) : (\mathbb{R}^n,0) \rightarrow (\mathbb{R}^k,0)$ be an analytic map, $2\leq k \leq n-1,$ $V=F^{-1}(0)$ and $\Sigma_F$ be the set of critical points of $F$, i.e., the set of points where the gradients $\nabla f_1,\ldots,\nabla f_k$ are linearly dependent. Of course, here and in the rest of the paper we assume that $F$ is not constant.

Let $\rho$ be the square of the distance function to the origin and denote by $\Sigma_{F,\rho}$ the set of critical points of the pair $(F,\rho)$, i.e., the set of points where the gradients $\nabla \rho,\nabla f_1,\ldots,\nabla f_k$ are linearly dependent.

\vspace{0.2cm}

It follows by definition that $\Sigma_F\subseteq \Sigma_{F,\rho}.$

\begin{definition}\cite{Ma}\label{defMilnorCond}Given $F$ and $\rho $ as above.
\vspace{0.2cm}

\begin{enumerate}
\item We say that $F$ satisfies Milnor's condition $(a)$ at the origin, if $\Sigma_F \subset V$ in a neighborhood of the origin.

 \item We say that $F$ satisfies Milnor's condition (b) at the origin, if $0$  is isolated in $V \cap \overline{\Sigma_{F,\rho}\setminus V}$ in a neighborhood of the origin.
\end{enumerate}
\end{definition}

Next example shows that the Milnor condition $(a)$ is not enough to ensure the existence of Milnor's fibration.

\begin{example} This example is inspired by examples in \cite{CSS}.

Let $f:(\mathbb{R}^{3},0)\to (\mathbb{R}^{2},0),$ $f(x,y,z)=(x^2z+y^2,x).$ It is easy to see that $\Sigma_f=\{(0,0,z): z\in \mathbb{R}\}\subseteq V$ and so Milnor's condition $(a)$ holds. However, for any $\delta > 0$ we have that the fibers $f^{-1}(\delta , 0)\neq f^{-1}(-\delta,0).$

\end{example}

\begin{remark}
It follows from definition $2)$ the equivalence:

The mapping $F$ satisfies Milnor's condition $(b)$ at origin if and only if there exist $\epsilon_{0}>0$ such that, for each
 $0<\epsilon \leq \epsilon_{0}$ we have $B_{\epsilon}\cap V \cap (\overline{\Sigma_{F,\rho}\setminus V}) \subseteq \{0\}$ if and only if
for each $\epsilon>0$ small enough, there exist $\delta>0$, $0<\delta \ll \epsilon $ such that the restriction map $F_|:S_{\epsilon}^{n-1}\cap F^{-1}({B^{p}_{\delta}\setminus \{0\}})\to B^{p}_{\delta}\setminus \{0\}$ is a smooth submersion (and onto, if the link of $F^{-1}(0)$ is not empty).
\end{remark}

We say that $\epsilon >0$ is a {\it Milnor radius for $F$ at origin}, provided that $B_{\epsilon}\cap {(\overline{\Sigma_{F}-V}})= \varnothing $, and $B_{\epsilon}\cap {V\cap(\overline{\Sigma_{F,\rho}\setminus V})}\subseteq \{0\}$, where
$B_{\epsilon}$ denotes the closed ball in $\mathbb{R}^{n}$ with radius $\epsilon$.

\vspace{0.2cm}

Consequently under Milnor's conditions $(a)$ and $(b)$, we can conclude that for all regular values close to the origin the respective fibers into the closed $\epsilon-$ball are smooth and transverse to the sphere $S^{n-1}_{\epsilon}.$

\begin{theorem}[\cite{Ma}, page 284, Theorem 4.3 ]

Let $F:(f_1,\ldots,f_k) : (\mathbb{R}^n,0) \rightarrow (\mathbb{R}^k,0)$ and $\epsilon_{0}>0$ be a Milnor's radius for $F$ at origin. Then, for each $0<\epsilon \leq \epsilon_{0}$, there exist $\delta,$ $0<\delta \ll \epsilon ,$ such that

\begin{equation}
F_{|}:B_{\epsilon}\cap F^{-1}(B^{p}_{\delta}\setminus \{0\})\to B^{p}_{\delta}\setminus \{0\}
\end{equation}
is the projection of a smooth locally trivial fiber bundle.

\end{theorem}

\proof (Idea)

Since $\epsilon_{0}>0$ is a Milnor's radius for $F$ at origin, we have that $\Sigma_{F}\cap B_{\epsilon_{0}} \subset V\cap B_{\epsilon_{0}}$.
It means that, for all $0<\epsilon \leq \epsilon_{0}$ the map $F_|: \mathring{B_{\epsilon }}\setminus V \to R^{k}$ is a smooth submersion in the open ball $\mathring{B_{\epsilon }}$.

From the Milnor condition $(b)$, and the remark above, it follows that: for each $\epsilon$ there exists $\delta $, $0<\delta \ll \epsilon $, such that
\begin{equation*}
 \displaystyle{ F_|: S_{\epsilon}^{n-1} \cap F^{-1}(B_{\delta}-\{0\}) \to B_{\delta}-\{0\}}
\end{equation*}
is a submersion on the boundary $\displaystyle{S_{\epsilon}^{n-1}}$ of the closed ball $\displaystyle{ B_{\epsilon}}$.

Now, combining these two conditions we have that, for each
$\epsilon $, we can choose $\delta $ such that

\begin{equation*}
 F_|: B_{\epsilon}\cap F^{-1}(B_{\delta}-\{0\}) \to  B_{\delta}-\{0\}
\end{equation*}
is a proper smooth submersion. Applying the Ehresmann Fibration Theorem for the manifold with boundary $B_{\epsilon}$, we get that it is a smooth locally trivial fibration. \endproof

\begin{corollary}
Let $F:(f_1,\ldots,f_k) : (\mathbb{R}^n,0) \rightarrow (\mathbb{R}^k,0)$ and $\epsilon_{0}>0$ be a Milnor's radius for $F$ at origin. Then, for each $0<\epsilon \leq \epsilon_{0}$, there exists $\delta,$ $0<\delta \ll \epsilon ,$ such that

\begin{equation}
F_{|}:B_{\epsilon}\cap F^{-1}(S^{p-1}_{\delta})\to S^{p-1}_{\delta}
\end{equation}
is the projection of a smooth locally trivial fiber bundle.

\end{corollary}

\begin{example}\label{R1}

Let $f:(\mathbb{C}^{n},0)\to (\mathbb{C},0)$ be a holomorphic function germ, then it satisfies the Milnor conditions $(a)$ and $(b).$

In fact, it can be seen as an application of \L ojasiewicz's inequality (see \cite{Lo}) which states that, in a small neighborhood of the origin, there are constants $C>0$ and $0<\theta<1$ such that

$$|f(x)|^{\theta}\leq C\|\nabla f(x)\|.$$

\vspace{0.2cm}

It is easy to see that Milnor condition $(a)$ holds. In \cite{HL}, page 323, Hamm and L\^e proved that the \L ojasiewicz inequality implies Thom $a_{f}-$condition for a Whitney $(a)$ stratification of $V.$ Therefore, Milnor's condition $(b)$ follows.

\end{example}

\begin{example}

Let $F:(\mathbb{R}^{n},0)\to (\mathbb{R}^{k},0)$ be an analytic map-germ with an isolated singular point at origin. Then, Milnor's conditions $(a)$ and $(b)$ above hold. In fact, Milnor's condition $(b)$ follows since the zero locus is transversal to all small spheres.
\end{example}

\section{Some results about subanalytic sets}

Let us recall some terminology and results  on the critical points of a function on the link of a real subanalytic set.
The situation is described as follows. Let $Y \subset \mathbb{R}^n$ be a smooth subanalytic set of dimension $d$ that contains $0$ in its closure and let $g : \mathbb{R}^n \rightarrow \mathbb{R}$ be a smooth subanalytic function such that $g(0)=0$.

\begin{lemma}\label{PtsCrit}
The critical points of $g_{\vert Y}$ lie in $\{g=0\}$ in a neighborhood of the origin.
\end{lemma}

 \proof By the Curve Selection Lemma, we can assume that there is a smooth subanalytic curve $p: [0,\nu[ \rightarrow \overline{Y}$ such that $p(0)=0$ and $p(t)$ is a critical point of $g_{\vert Y}$ for $t\in ]0,\nu[$.
Therefore we have
$$(g \circ p(t))' =\langle \nabla g(p(t)), p'(t) \rangle =0,$$
since $p'(t)$ is a tangent vector to $Y$ at $p(t)$. This implies that $g \circ p(t)=g(p(0))=0$. \endproof

Now we are interested in the critical points of $g_{\vert Y \cap S_\varepsilon}$, where $0 < \varepsilon \ll 1$,  lying in $\{g \not= 0 \}$. Let $q$ be such a critical point.
By the previous lemma, we know that $\nabla g_{\vert Y}(q) \not= 0$ and so there exists $\lambda (q) \not= 0$ such that
$$\nabla g_{\vert Y} (q) = \lambda (q) \nabla \rho_{\vert Y} (q).$$
\begin{definition}\label{outin}
We say that $q \in \{Y \cap S_ \varepsilon \}$ is an outwards-pointing (resp. inwards-pointing) critical point for $g_{\vert Y \cap S_\varepsilon}$ if $\lambda(q)>0$ (resp. $\lambda(q)<0$).
\end{definition}

\begin{lemma}
The point $q \in \{g \not= 0 \}$ is an outwards-pointing (resp. inwards-pointing) critical point for $g_{\vert Y \cap S_\varepsilon}$ if and only if $g(q)> 0$ (resp. $g(q)<0$).
\end{lemma}

\proof Let us assume that $\lambda (q) >0$. By the Curve Selection Lemma, there exists a smooth subanalytic curve $p:[0,\nu[ \rightarrow \overline{Y}$ passing through $q$ such that $p(0)=0$
and for $t\not= 0$, $p(t)$ is a critical point of $g_{\vert Y \cap S_{\Vert p(t) \Vert}}$ with $\lambda(p(t)) >0$. Therefore we have:
$$(g \circ p)'(t) =\langle \nabla g _{\vert Y}(p(t)), p'(t) \rangle= \lambda(p(t)) \langle \nabla \rho _{\vert Y}(p(t)), p'(t) \rangle= \lambda(p(t)) (\rho  \circ p)'(t).$$
But $(\rho  \circ p)'>0$ for otherwise $(\rho  \circ p)' \le 0$ and $\rho \circ p$ would be decreasing.  Since $\rho (p(t))$ tends to $0$ as $t$ tends to $0$, this would imply that $\rho  \circ p(t) \le 0$, which is impossible. Hence we can conclude that $(g \circ p)'>0$ and $g \circ p$ is strictly increasing. Since $g \circ p(t)$ tends to $0$ as $t$ tends to $0$, we see that $g \circ p (t) >0$ for $t>0$. Similarly if $\lambda (q) <0$ then $g(q)<0$.  $\hfill \Box$

\section{Basic results on Milnor's conditions $(a)$ and $(b)$}

In this section we consider $F=(f_1,\ldots,f_k) : (\mathbb{R}^n,0) \rightarrow (\mathbb{R}^k,0),$ $1\leq k \leq n-1,$ an analytic mapping.
Let us consider $l \in \{1,\ldots, k\}$ and $I= \{i_1,\ldots,i_l\}$ an $l$-tuple of pairwise distinct elements of $ \{1,\ldots,k\}$ and let us denote by  $f_I$ the mapping $(f_{i_1},\ldots,f_{i_l}) :
(\mathbb{R}^n,0) \rightarrow (\mathbb{R}^l,0)$. Suppose that $F$ satisfies Milnor condition $(a)$ at the origin. Then, we have
$$\Sigma_{f_I} \subset \Sigma_F \subset F^{-1}(0) \subset f_{I}^{-1}(0),$$
and so by definition the map $f_I$ also satisfies Milnor's condition (a) at the origin.

\vspace{0.3cm}

It is clear that $\Sigma_{f_I,\rho} \subset \Sigma_{F,\rho}.$ We will show below that if $F$ satisfies Milnor's condition $(b)$ at the origin, then any mapping $f_I$ also satisfies Milnor's condition $(b)$ at the origin.

\begin{lemma}\label{lemmaCondb} Assume that $F$ satisfies Milnor's conditions $(a)$ and $(b)$ at the origin. Then, for $l \in \{1,\ldots, k\}$ and $I= \{i_1,\ldots,i_l\} \subset \{1,\ldots,k\}$, the maps $f_I :  (\mathbb{R}^n,0) \rightarrow (\mathbb{R}^l,0)$ satisfies Milnor's conditions $(a)$ and $(b)$.
\end{lemma}
\proof If a map $f_I$ does not satisfy condition (b), then $0$ is not isolated in $f_I^{-1}(0) \cap \overline{\Sigma_{f_I,\rho} \setminus f_I^{-1}(0)}$. This implies that there exists a sequence of points $(y_n)_{n \in \mathbb{N}}$ tending to the origin such that $f_I(y_n)=0$ and $y_n$ belongs to $\overline{\Sigma_{f_I,\rho} \setminus f_I^{-1}(0)}$.

If $y_n$ belongs to $\overline{\Sigma_{f_I,\rho} \setminus f_I^{-1}(0)}$, there exists a sequence of points $(y_n^k)_{k \in \mathbb{N}}$  tending to $y_n$ such that $f_I(y_n^k) \not= 0$ and the gradients $\nabla \rho,\nabla f_{i_1}, \ldots,\nabla f_{i_l}$ are linearly dependent at the points $y_n^k$. Hence the gradients $\nabla \rho,\nabla f_{i_1}, \ldots,\nabla f_{i_l}$ are also linearly dependent at $y_n$.
But, since in a neighborhood of the origin $\rho$ has no critical point on $f_I^{-1}(0) \setminus \Sigma_{f_I}$ by Lemma 3.1, we see that $y_n$ belongs to $\Sigma_{f_I}$ if $n$ is big enough and so, $y_n$ belongs to $V$ because
$\Sigma_{f_I} \subset \Sigma_F \subset V$.

On the other hand, the points $y_n^k$'s belong to $\Sigma_{F,\rho} \setminus V$ as well, so $y_n$ lies in $\overline{\Sigma_{F,\rho} \setminus V }$ and therefore $0$ is not isolated in  $V \cap \overline{\Sigma_{F,\rho} \setminus V}$.
\endproof

\begin{corollary}
There exists $\epsilon_{0}>0$ such that, for all $l \in \{1,\ldots, k\}$ and $I= \{i_1,\ldots,i_l\} \subset \{1,\ldots,k\}$, the maps $f_I :  (\mathbb{R}^n,0) \rightarrow (\mathbb{R}^l,0)$ have $\epsilon_{0}$ as a Milnor's radius.

\end{corollary}

\section{Critical points on the boundary of the Milnor fibre}
From now on, we consider an analytic mapping $F:(f_1,\ldots,f_k) : (\mathbb{R}^n,0) \rightarrow (\mathbb{R}^k,0)$ with a Milnor's radius $\epsilon >0$, $V=F^{-1}(0)$, the mapping $\phi=(f_1,\ldots,f_{k-1})$ and $g=f_k$. By Lemma \ref{lemmaCondb}, these two maps also satisfy Milnor's conditions (a) and (b).

\vspace{0.2cm}

In this section we will study the behaviour of the critical points of the function $g$ restricted to the boundary of the Milnor fibre of the mapping $\phi$.

The next lemma is inspired by \cite{Sz1}, pages 411--412.

\begin{lemma}\label{LemSigma}
There exist a positive constant $C$ and an integer $N$ such that $\Vert F(x) \Vert \ge C \Vert x \Vert^N$ for every $x \in \Sigma_{F,\rho} \setminus V$ sufficiently close to the origin.
\end{lemma}

\proof Let
$$\Gamma = \left\{ (x,r,y) \in \mathbb{R}^n \times \mathbb{R} \times \mathbb{R}^k \ \vert \ \rho(x)=r, x \in \Sigma_{F,\rho} \hbox{ and } y=F(x) \right\},$$
and let $\pi :  \mathbb{R}^n \times \mathbb{R} \times \mathbb{R}^{k} \rightarrow   \mathbb{R} \times \mathbb{R}^k$ be the projection on the last  $k+1$ components.
Since $\pi_{\vert \Gamma}$ is proper, then $Z=\pi(\Gamma)$ is a closed semi-analytic set. Let us write $Z_1= \mathbb{R} \times \{0\} \subset \mathbb{R} \times \mathbb{R}^k$ and let $Z_2$ be the closure of $Z \setminus Z_1$.
Then $0$ is isolated in $Z_1 \cap Z_2$. If it is not the case, this means that there is a sequence of points $z_i=(r_i,0)$ in $Z_1$ tending to $0$ such that  $z_i$ belongs to $Z_2$.
Hence for each $i$, there is a sequence of points $(z_i^j)_{j \in \mathbb{N}}$ in $Z \setminus Z_1$ tending to $z_i$. Let us write $z_i^j=(r_i^j,y_i^j)$. Since $z_i^j$ is in $Z \setminus Z_1$, this implies that there exists $x_i^j$ in $\Sigma_{F,\rho}$ such that $F(x_i^j)=y_i^j$. Taking a subsequence if necessary, we can assume that $(x_i^j)$ tends to a point $x_i$ which belongs to $V$ because $F(x_i^j)=y_i^j$ tends to $0$ and such that $\rho(x_i)=r_i$ because $r_i^j=\rho(x_i^j)$ tends to $r_i$. But since $r_i$ tends to $0$, this implies that $0$ is not isolated in $V \cap \overline{\Sigma_{F,\rho} \setminus V}$, which contradicts Milnor's condition (b).

By the \L ojasiewicz inequality, there exist a constant $C >0$ and an integer $N >0$ such that
$$ \Vert y \Vert  \ge C r ^N,$$
for $(r,y) \in Z_2$ sufficiently close to the origin. So if $x \in \Sigma_{F,\rho}$ and $F(x) \not= 0$, then $\Vert F(x) \Vert \ge C \rho (x) ^N$ if $\Vert x \Vert$ is small enough. \endproof

As a consequence, we see that for $\epsilon >0$ small enough, there exists $\delta_\epsilon >0$ such that if $0 < \Vert \delta \Vert < \delta_\epsilon$, then $F^{-1}(\delta)$ intersects $S_\epsilon$ transversally.

\begin{remark}
The lemma above can the generalized by changing the square of the distance function to the origin with any subanalytic function $\rho$ smooth, positive and proper, such that locally $\rho^{-1}(0)=0.$

\end{remark}

\begin{corollary}\label{CorCriticPtsBound}
For $\epsilon >0$ small enough, there exists  $\delta_\epsilon >0$ such that if $0 < \Vert \delta \Vert < \delta_\epsilon$, then the critical points of $g_{\vert \phi^{-1}(\delta) \cap S_\epsilon}$ lie in
$\{ \vert g \vert \ge \frac{\sqrt{3}}{2}C \epsilon^N \}$.
\end{corollary}

 \proof Applying the above lemma and consequence to the mapping $\phi=(f_1,\ldots,f_{k-1})$, we see that there exist a constant $D>0$ and an integer $M >0$ such that
$$\Vert \phi (x) \Vert \ge D \Vert x \Vert^M,$$
for $x \in \Sigma_{\phi,\rho} \setminus \phi^{-1}(0)$ sufficiently close to the origin.
Let us fix $\epsilon>0$ sufficiently small so that $S_\epsilon$ intersects $\phi^{-1}(0) \setminus \Sigma_\phi$ transversally. If $\delta  \in \mathbb{R}^{k-1}$ is such that $0< \Vert \delta \Vert \le \frac{D}{2} \epsilon^M$, then
$S_\epsilon$ intersects the fibre $\phi^{-1}(\delta)$ transversally by the above inequality.
Let us choose $\delta \in \mathbb{R}^{k-1}$ such that $0< \Vert \delta \Vert \le \hbox{Min}\{ \frac{D}{2} \epsilon^M, \frac{C}{2} \epsilon^N \}$. If $x$ is a critical point $g_{\vert \phi^{-1}(\delta) \cap S_\epsilon}$ then
$x$ belongs to $\Sigma_{F,\rho} \setminus V$ and so $\Vert F(x) \Vert \ge C \epsilon^N$. This implies that
$$g(x)^2 \ge C^2 \epsilon^{2N} - \Vert \phi \Vert^2 \ge C^2\epsilon^{2N} -\frac{C^2}{4} \epsilon^{2N},$$
and so $\vert g(x) \vert \ge \frac{\sqrt{3}}{2}C \epsilon^{N}$. \endproof

Now we can return to our map $F:(\phi,g) : (\mathbb{R}^n,0)\rightarrow (\mathbb{R}^k,0)$. Let us choose $\epsilon$ and $\delta \in \mathbb{R}^{k-1}$ such that $0 < \Vert \delta \Vert \ll \epsilon \ll 1$ and the critical points of $g_{\vert \phi^{-1}(\delta) \cap S_\epsilon}$ lie in $\{ \vert g \vert \ge \frac{\sqrt{3}}{2}C \epsilon^N \}$. \endproof

\begin{lemma}
The critical points of $g_{\vert \phi^{-1}(\delta) \cap S_\epsilon}$  in $\{g \ge \frac{\sqrt{3}}{2}C \epsilon^N \}$ are outwards-pointing and the critical points of $g_{\vert \phi^{-1}(\delta) \cap S_\epsilon}$  in $\{g \le
- \frac{\sqrt{3}}{2}C \epsilon^N\}$ are inwards-pointing.
\end{lemma}

\proof Let us prove the statement about the critical points in $\{g \ge \frac{\sqrt{3}}{2}C \epsilon^N \}$.
 Let us remark first that such a critical point is not a critical point of $g_{\vert \phi^{-1}(\delta)}$ because $\Sigma_{\phi,g} \subset \phi^{-1}(0) \cap g^{-1}(0)$ in a neighborhood of the origin.
Therefore if the statement is not verified, then this means that we can find a sequence of points $(q_n)_{n \in \mathbb{N}}$ in $S_\epsilon \cap \{g \ge \frac{\sqrt{3}}{2}C \epsilon^N \}$ such that $\phi(q_n)$ tends to $0$ and $q_n$ is an inwards-pointing critical point of $g_{\vert \phi^{-1}(\phi(q_n)) \cap S_\epsilon}$. Taking a subsequence if necessary, this produces a point $q$ in $\phi^{-1}(0) \cap S_\epsilon$ such that $g(q) \ge  \frac{\sqrt{3}}{2}C \epsilon^N$ and $q$ is a critical point of $g$ restricted to $\phi ^{-1}(0) \setminus \Sigma_{\phi} \cap S_\epsilon$, because
$$\Sigma_{\phi} \subset \Sigma_F \subset F^{-1}(0) \subset g^{-1}(0).$$
Futhermore as the limit of a sequence of inwards-pointing critical points, $q$ is either inwards-pointing, which is impossible by the previous lemma, or is a critical point of $g$ restricted to $\phi^{-1}(0) \setminus \Sigma_\phi$, which is also impossible by Lemma \ref{PtsCrit}. \endproof

\section{Topology of the Milnor fibre}

We keep the notations of the previous section.
We denote by $M_F$ the Milnor fibre of the mapping $F$ and by $M_\phi$ the Milnor fibre of $\phi$. Here we need to assume that $k \ge 3$ so that the Milnor fibre of $\phi$ is well-defined.

\begin{theorem} [Milnor's conjecture, \cite{Mi}, page 100] \label{Milnorconjec}
The fibre  $M_\phi$ is homeomorphic to $M_F \times [-1,1] $.
\end{theorem}
{\it Proof}. Let us choose $\epsilon \in \mathbb{R}$ and $\delta \in \mathbb{R}^k$ such that:
\begin{enumerate}
\item $0< \Vert \delta \Vert \ll \epsilon \ll 1$,
\item $M_\phi$ is homeomorphic to $\phi^{-1}(\delta) \cap B_\epsilon$,
\item $M_F$ is homeomorphic to $\phi^{-1}(\delta) \cap g^{-1}(0) \cap B_\epsilon$,
\item the critical points of $g$ restricted to $\phi^{-1}(\delta) \cap S_\epsilon$ lie in $\{g \not= 0\}$, are outwards-pointing in $\{g > 0\}$ and inwards-pointing in $\{ g < 0 \}$.
\end{enumerate}
Note that $g_{\vert \phi^{-1}(\delta) \cap \mathring{B_\epsilon}}$ has no critical points because $\Sigma_{\phi,g} \subset \phi^{-1}(0) \cap g^{-1}(0)$.
Then by Morse theory for manifolds with boundary (see \cite{H}), $\phi^{-1}(\delta) \cap B_\epsilon \cap \{ g \ge 0 \}$ is homeomorphic to $\phi^{-1}(\delta) \cap B_\epsilon \cap \{ g = 0 \} \times [0,1]$ and $\phi^{-1}(\delta) \cap B_\epsilon \cap \{ g \le 0 \}$ is homeomorphic to $\phi^{-1}(\delta) \cap B_\epsilon \cap \{ g = 0 \} \times [-1,0]$.
Therefore $\phi^{-1}(\delta) \cap B_\epsilon$ is homeomorphic to  $\phi^{-1}(\delta) \cap B_\epsilon \cap \{ g = 0 \} \times [-1,1]$ because it is homeomorphic to the gluing of $\phi^{-1}(\delta) \cap B_\epsilon \cap \{ g \ge 0 \}$ and $\phi^{-1}(\delta) \cap B_\epsilon \cap \{ g \le 0 \}$ along $\phi^{-1}(\delta) \cap B_\epsilon \cap \{ g = 0 \}$. $\hfill \Box$

\begin{corollary}
Under the above conditions we have $\chi(M_F)=\chi(M_\phi)$.
\end{corollary}

Proceeding by induction on the number of components of the mapping, we can easily prove the following corollary.
\begin{corollary}\label{CharFib1}
Let $l \in \{2,\ldots,k\}$ and let $I=\{i_1,\ldots,i_l\}$ be an $l$-tuple of pairwise distinct elements of $ \{1,\ldots,k\}$. Then we have $\chi(M_{f_I})=\chi(M_F)$.
\end{corollary}

It remains to consider the fibres of the function $f_j : (\mathbb{R}^n,0) \rightarrow (\mathbb{R},0)$. Here such a function admits two Milnor fibres : $M^+_{f_{\{j\}}}= f_j^{-1}(\delta) \cap B_\epsilon$ and
$M^-_{f_{\{j\}}}= f_j^{-1}(-\delta) \cap B_\epsilon$, where $0 < \delta \ll \epsilon \ll 1$.

Let us write for instance $f=f_1$ and $g=f_2$.
Using the same argument as above, we see that $M^+_f$ is homeomorphic to $M_{(f,g)}\times [-1,1]$ and that $M^-_f$ is also homeomorphic $M_{(f,g)} \times [-1,1]$.

\begin{corollary}\label{CharFib2}
For every $j \in \{1,\ldots,k \}$, we have $\chi(M^+_{f_{\{j\}}})=\chi(M^-_{f_{\{j\}}})=\chi(M_F)$.
\end{corollary}
\endproof

\section{Topology of the  links}
In this section we give several relations between the Euler characteristics of the links of $f_I^{-1}(0)$ and the Euler characteristic of the Milnor fibre of $F$.

Let us choose $l \in \{1,\ldots,k \}$ and an  $l$-tuple $I=\{i_1,\ldots,i_l\}$ of pairwise distinct elements of $\{1,\ldots,k\}$. We write $J=\{i_1,\ldots,i_{l-1}\}$ and $g=f_{i_l}$.
We also denote by $\mathcal{L}_I$ (resp. $\mathcal{L}_J$) the link of the zero-set of $f_I$ (resp. $f_J$). If $l=1$ then $J= \emptyset$ and we put $f_J=0$.
\begin{proposition}\label{CharL1}
We have:
$$\chi(\mathcal{L}_J)-\chi(\mathcal{L}_I)=(-1)^{n-l} 2 \chi(M_F).$$
\end{proposition}
 \proof Let us write $V_J=f_J^{-1}(0)$. By a deformation argument due to Milnor, $V_J \cap g^{-1}(\delta) \cap B_\epsilon$ is homeomorphic to $V_J \cap \{ g \ge \delta \} \cap S_\epsilon$
and $V_J \cap g^{-1}(-\delta) \cap B_\epsilon$ is homeomorphic to $V_J \cap \{ g \le -\delta \} \cap S_\epsilon$ for $0 < \delta \ll \epsilon \ll 1$.
By the Mayer-Vietoris sequence, we can write:
$$\chi(V_J \cap S_\epsilon) = \chi (V_J \cap S_\epsilon \cap \{ g \ge \delta \} )+\chi (V_J \cap S_\epsilon \cap \{ g \le -\delta \} )+ $$
$$\chi (V_J \cap S_\epsilon \cap  \{ -\delta \le g \le \delta \} ) -\chi(V_J \cap S_\epsilon \cap  \{ g = \delta \} )-\chi(V_J \cap S_\epsilon \cap \{ g = -\delta \} ).$$
By the above remark and Corollaries \ref{CharFib1} and \ref{CharFib2}, the first two terms of the right-hand side of this equality are equal to $\chi(M_F)$.
The third term is equal to $\chi(\mathcal{L}_I)$ because by Durfee's result [Dur], $\mathcal{L}_I$ is a retract by deformation of $V_J \cap S_\epsilon \cap  \{ -\delta \le g \le \delta \} $.
Furthermore, if $n-l$ is even then the two last Euler characteristics are equal to $0$ because $V_J \cap S_\epsilon \cap  \{ g = \delta \} $ and $V_J \cap S_\epsilon \cap  \{ g =- \delta \} $ are odd-dimensional compact manifolds. If $n-l$ is odd, they are equal to $2\chi(M_F)$ because they are boundaries of odd-dimensional Milnor fibres of $f_I$. \endproof

\begin{corollary}
Let $j \in \{1,\ldots,k\}$. If $n$ is even, then we have $\chi(\mathcal{L}_{\{j\}})= 2 \chi(M_F)$ and if $n$ is odd, then we have $\chi(\mathcal{L}_{\{j\}})=2-2\chi(M_F)$.
\end{corollary}
\proof We apply the previous proposition to the case $l=1$. In this case, if $n$ is even then $\chi(\mathcal{L}_J)=0$ and if $n$ is odd then $\chi(\mathcal{L}_J)=2$. $\hfill \Box$

\begin{corollary}
Let $l\in \{3,\ldots,k \}$ and let $I=\{i_1,\ldots,i_l\} \subset \{1,\ldots,k\}$. Let $K$ be an $(l-2)$-tuple of pairwise distincts elements of $I$.
Then we have: $\chi(\mathcal{L}_K)= \chi(\mathcal{L}_J)$.
\end{corollary}
\proof Let $J$ be an $(l-1)$-tuple built form adding to $K$ one element of $I \setminus K$. By the previous proposition, we see that $\chi(\mathcal{L}_J)-\chi(\mathcal{L}_I)=\chi(\mathcal{L}_J)-\chi(\mathcal{L}_K)$.
\endproof

So, in order to express the Euler characteristics of all the links $\mathcal{L}_I$, we just need to compute the Euler characteristic of a link $\mathcal{L}_I$ where $\# I=2$.
Let us set $I=\{1,2\}$. By Proposition \ref{CharL1}, we find that
$$\chi(\mathcal{L}_I)=\chi(\mathcal{L}_{\{1\}})-(-1)^n 2\chi(M_F).$$
So if $n$ is even, we see that $\chi(\mathcal{L}_I)=0$ and if $n$ is odd,  we see that $\chi(\mathcal{L}_I)=2$.
We can summarize all these results in the following theorem.

\begin{theorem}\label{CharLink}
Let $l\in \{1,\ldots,k \}$ and let $I=\{i_1,\ldots,i_l\}$ be an $l$-tuple of pairwise distinct elements of $\{1,\ldots,k\}$.
If $n$ is even, then we have:
$$\chi(\mathcal{L}_I)= 2\chi(M_F) \hbox{ if } l \hbox{ is odd and }\chi(\mathcal{L}_I)=0 \hbox{ if } l \hbox{ is even}.$$
If $n$ is odd, then we have:
$$\chi(\mathcal{L}_I)=2- 2\chi(M_F) \hbox{ if } l \hbox{ is odd and }\chi(\mathcal{L}_I)=2 \hbox{ if } l \hbox{ is even}.$$
\end{theorem}
$\hfill \Box$

\section{Topology of related semi-analytic sets}
In this section, we establish formulas for the Euler characteristics of several semi-analytic sets defined from the components of the map $F$.
We are interested first in the sets of the form
$$f_I ^{-1}(\delta) \cap \{f_{j_1} \epsilon_1 0, \ldots, f_{j_s} \epsilon_s 0 \} \cap B_\epsilon,$$
where $I=\{i_1,\ldots,i_l \}$ is an $l$-tuple of pairwise distinct elements of $\{1,\ldots,k \}$, $l+s \le k$,  $\delta$ is a sufficiently small regular value of $f_I$, $j_1,\ldots,j_s$ are pairwise distinct elements of $\{1,\ldots,k\} \setminus I$ and for $i \in \{1,\ldots, s \}$, $\epsilon_i \in \{\le,\ge \}$.

\begin{proposition}\label{CharSemi1}
We have:
$$\chi \left( f_I ^{-1}(\delta) \cap \{f_{j_1} \epsilon_1 0, \ldots, f_{j_s} \epsilon_s 0 \} \cap B_\epsilon \right) = \chi (M_F).$$
\end{proposition}
\proof We prove the result by induction on $s$. Let us state the induction hypothesis IH($s$) properly.

\vskip0,2cm

\noindent Let $s \in \{1,\ldots,k \}$. For any $l \in \{1, \ldots, k \}$ such that $l+s \le k$, for any $l$-tuple $I=\{i_1,\ldots,i_l\}$ of pairwise distinct elements of $\{1,\ldots,k \}$, for any $s$-tuple $\{j_1,\ldots,j_s\}$ of pairwise distinct elements of
$\{1,\ldots,k \} \setminus I$, we have
$$\chi \left( f_I ^{-1}(\delta) \cap \{f_{j_1} \epsilon_1 0, \ldots, f_{j_s} \epsilon_s 0 \} \cap B_\epsilon \right) = \chi (M_F),$$
where $\delta$ is a sufficiently small regular value of $f_I$ and for $i \in \{1,\ldots, s \}$, $\epsilon_i \in \{\le,\ge \}$.

\vskip0,2cm

Let us prove IH(1).
By the same argument of Theorem \ref{Milnorconjec}, we can apply Morse theory for manifolds with boundary to $f_{j_1}$ restricted to $f_I^{-1}(\delta) \cap  B_\epsilon$, to find that
$$ \chi \left( f_I^{-1}(\delta) \cap  \{ f_{j_1} \ge 0 \} \cap B_\epsilon \right) -
 \chi \left( f_I^{-1}(\delta) \cap   \{ f_{j_1} = 0 \} \cap B_\epsilon \right) =0$$ and
$$\chi \left( f_I^{-1}(\delta) \cap  \{ f_{j_1} \le 0 \} \cap B_\epsilon \right) -
\chi \left( f_I^{-1}(\delta) \cap  \{ f_{j_1} = 0 \} \cap B_\epsilon \right) =0. $$
Now by Corollary \ref{CharFib1} and Corollary \ref{CharFib2} we get the result.

\vspace{0.2cm}

Let us assume that IH($s-1$) is satisfied and let us prove IH($s$).
By Morse theory for manifold with corners (see \cite{Du3}) applied to $f_{j_s}$ restricted to $f_I^{-1}(\delta) \cap  \{f_{j_1} \epsilon_1 0, \ldots, f_{j_{s-1}} \epsilon_{s-1} 0 \} \cap B_\epsilon$, we find that
$$\displaylines{
\qquad \chi \left( f_I^{-1}(\delta) \cap  \{f_{j_1} \epsilon_1 0, \ldots, f_{j_{s-1}} \epsilon_{s-1} 0, f_{j_s} \ge 0 \} \cap B_\epsilon \right) - \hfill \cr
\hfill \chi \left( f_I^{-1}(\delta) \cap  \{f_{j_1} \epsilon_1 0, \ldots, f_{j_{s-1}} \epsilon_{s-1} 0, f_{j_s} = 0 \} \cap B_\epsilon \right) =0, \qquad \cr}$$ and
$$\displaylines{
\qquad \chi \left( f_I^{-1}(\delta) \cap  \{f_{j_1} \epsilon_1 0, \ldots, f_{j_{s-1}} \epsilon_{s-1} 0, f_{j_s} \le 0 \} \cap B_\epsilon \right) - \hfill \cr
\hfill \chi \left( f_I^{-1}(\delta) \cap  \{f_{j_1} \epsilon_1 0, \ldots, f_{j_{s-1}} \epsilon_{s-1} 0, f_{j_s} = 0 \} \cap B_\epsilon \right) =0. \qquad \cr}$$
But, by the induction hypothesis IH($s-1$), we know that
$$\chi \left( f_I^{-1}(\delta) \cap  \{f_{j_1} \epsilon_1 0, \ldots, f_{j_{s-1}} \epsilon_{s-1} 0, f_{j_s} = 0 \} \cap B_\epsilon \right) = \chi(M_F).$$
\endproof

Now we look at the sets of the form
$$f_I^{-1}(0)  \cap \{f_{j_1} \epsilon_1 0, \ldots, f_{j_{s-1}} \epsilon_{s-1} 0, f_{j_s} = \delta \} \cap S_\epsilon,$$
where $\delta$ is a sufficiently small regular value of $f_{j_s}$.

\begin{lemma}
We have:
$$\chi \left( f_I^{-1}(0)  \cap \{f_{j_1} \epsilon_1 0, \ldots, f_{j_{s-1}} \epsilon_{s-1} 0, f_{j_s} = \delta \} \cap S_\epsilon\right) = $$
$$\chi \left( f_I^{-1}(0)  \cap \{f_{j_1} \epsilon_1 0, \ldots, f_{j_{s-1}} \epsilon_{s-1} 0, f_{j_s} = -\delta \} \cap S_\epsilon\right) .$$
\end{lemma}
\proof We prove the result by induction on $s$. Let us state the induction hypothesis IH($s$) properly.

\vskip0,2cm

\noindent Let $s \in \{1,\ldots,k \}$. For any $l \in \{1, \ldots, k \}$ such that $l+s \le k$, for any $l$-tuple $I=\{i_1,\ldots,i_l\}$ of pairwise distinct elements of $\{1,\ldots,k \}$, for any $s$-tuple $\{j_1,\ldots,j_s\}$ of pairwise distinct elements of
$\{1,\ldots,k \} \setminus I$, we have
$$\chi \left( f_I^{-1}(0)  \cap \{f_{j_1} \epsilon_1 0, \ldots, f_{j_{s-1}} \epsilon_{s-1} 0, f_{j_s} = \delta \} \cap S_\epsilon\right) = $$
$$\chi \left( f_I^{-1}(0)  \cap \{f_{j_1} \epsilon_1 0, \ldots, f_{j_{s-1}} \epsilon_{s-1} 0, f_{j_s} = -\delta \} \cap S_\epsilon\right),$$
where $\delta$ is a sufficiently small regular value of $f_{j_s}$ and for $i \in \{1,\ldots, s-1 \}$, $\epsilon_i \in \{\le,\ge \}$.

\vskip0,2cm

Let us prove first IH(1). This is easy because $f_I^{-1}(0) \cap \{f_{j_1}=\delta\} \cap S_\epsilon$ is the boundary of the manifold $f_I^{-1}(0) \cap \{f_{j_1}=\delta\} \cap B_\epsilon$ and so its Euler characteristic is $0$ if dim$f_I^{-1}(0)$ is odd and it is $2 \chi(M_F)$ if dim$f_I^{-1}(0)$ is even.

Let us assume that IH($1$), $\ldots$, IH($s-1$) are satisfied and let us prove IH($s$).  If dim$f_I^{-1}(0)$ is odd then dim$f_I^{-1}(0) \cap \{f_{j_1} \epsilon_1 0, \ldots, f_{j_{s-1}} \epsilon_{s-1} 0, f_{j_s} = \delta \} \cap S_\epsilon $ is also odd.
But $f_I^{-1}(0) \cap \{f_{j_1} \epsilon_1 0, \ldots, f_{j_{s-1}} \epsilon_{s-1} 0, f_{j_s} = \delta \} \cap S_\epsilon $ is an odd-dimensional manifold with corners so, after rounding the corners, we can write
$$\chi \left( f_I^{-1}(0) \cap \{f_{j_1} \epsilon_1 0, \ldots, f_{j_{s-1}} \epsilon_{s-1} 0, f_{j_s} = \delta \} \cap S_\epsilon \right)=$$
$$\frac{1}{2} \chi  \Big(  \partial \left(f_I^{-1}(0) \cap \{f_{j_1} \epsilon_1 0, \ldots, f_{j_{s-1}} \epsilon_{s-1} 0, f_{j_s} = \delta \} \cap S_\epsilon \right) \Big).$$
By the Mayer-Vietoris sequence, we see that the Euler characteristic of the right-hand side is a linear combination, with coefficients equal to $\pm 1$, of the Euler characteristics of sets of the form:
$$f_I^{-1}(0) \cap \{f_{j_1} \nu_1 0, \ldots, f_{j_{s-1}} \nu_{s-1} 0, f_{j_s} = \delta \} \cap S_\epsilon,$$
where $\nu_i  \in \{\le,=,\ge \}$ for $i=1,\ldots,s-1$ and at least one of the $\nu_i$'s is the sign $=$. Hence we can apply the induction hypothesis to obtain the result.

If dim$f_I^{-1}(0)$ is even then dim$f_I^{-1}(0) \cap \{f_{j_1} \epsilon_1 0, \ldots, f_{j_{s-1}} \epsilon_{s-1} 0, f_{j_s} = \delta \} \cap S_\epsilon $ is also even.
But $f_I^{-1}(0) \cap \{f_{j_1} \epsilon_1 0, \ldots, f_{j_{s-1}} \epsilon_{s-1} 0, f_{j_s} = \delta \} \cap B_\epsilon $ is an odd-dimensional manifold with corners so, after rounding the corners, we can write
$$\chi \left( f_I^{-1}(0) \cap \{f_{j_1} \epsilon_1 0, \ldots, f_{j_{s-1}} \epsilon_{s-1} 0, f_{j_s} = \delta \} \cap B_\epsilon \right)=$$
$$\frac{1}{2} \chi  \Big(  \partial \left(f_I^{-1}(0) \cap \{f_{j_1} \epsilon_1 0, \ldots, f_{j_{s-1}} \epsilon_{s-1} 0, f_{j_s} = \delta \} \cap B_\epsilon \right) \Big).$$
We know by the previous proposition that the left-hand side of this equality does not depend on the sign of $\delta$. Let us examine the Euler characteristic of the right-hand side. By the Mayer-Vietoris sequence, it is equal to a linear combination, with coefficients equal to $\pm 1$, of the Euler characteristic of
$$f_I^{-1}(0) \cap \{f_{j_1} \epsilon_1 0, \ldots, f_{j_{s-1}} \epsilon_{s-1} 0, f_{j_s} = \delta \} \cap S_\epsilon ,$$
the Euler characteristics of sets of the type:
$$f_I^{-1}(0) \cap \{f_{j_1} \nu_1 0, \ldots, f_{j_{s-1}} \nu_{s-1} 0, f_{j_s} = \delta \} \cap B_\epsilon,$$
where $\nu_i  \in \{\le,=,\ge \}$ for $i=1,\ldots,s-1$,
and the Euler characteristics of sets of the type:
$$f_I^{-1}(0) \cap \{f_{j_1} \nu_1 0, \ldots, f_{j_{s-1}} \nu_{s-1} 0, f_{j_s} = \delta \} \cap S_\epsilon,$$
where $\nu_i  \in \{\le,=,\ge \}$ for $i=1,\ldots,s-1$ and at least one of the $\nu_i$'s is the sign $=$. Since the Euler characteristics of the sets of the first type do not depend on the sign of $\delta$ by the previous proposition, and those of the sets of the second type do not neither by the induction hypothesis, we get the result.  $\hfill \Box$

Now we study the links of the sets of the form
$$f_I ^{-1}(0) \cap \{f_{j_1} \epsilon_1 0, \ldots, f_{j_s} \epsilon_s 0 \} ,$$
where $I=\{i_1,\ldots,i_l \} \subset \{1,\ldots,k \}$, $l+s \le k$,  $j_1,\ldots,j_s$ are pairwise distinct elements of $\{1,\ldots,k\} \setminus I$ and for $i \in \{1,\ldots, s \}$, $\epsilon_i \in \{\le,\ge \}$.
\begin{theorem}
We have:
$$\chi \left( \mathcal{L}_{I} \cap \{f_{j_1} \epsilon_1 0, \ldots, f_{j_s} \epsilon_s 0 \} \right)= \chi(M_F),$$
if $n$ is even and
$$\chi \left( \mathcal{L}_{I} \cap \{f_{j_1} \epsilon_1 0, \ldots, f_{j_s} \epsilon_s 0 \} \right) =2- \chi(M_F),$$
if $n$ is odd.
\end{theorem}
\proof  We prove the result by induction on $s$. Let us state the induction hypothesis IH($s$) properly.

\vskip0,2cm

\noindent Let $s \in \{1,\ldots,k \}$. For any $l \in \{1, \ldots, k \}$ such that $l+s \le k$, for any $l$-tuple $I=\{i_1,\ldots,i_l\}$ of pairwise distinct elements of $\{1,\ldots,k \}$, for any $s$-tuple $\{j_1,\ldots,j_s\}$ of pairwise distinct elements of
$\{1,\ldots,k \} \setminus I$, we have
$$\chi \left( \mathcal{L}_{I} \cap \{f_{j_1} \epsilon_1 0, \ldots, f_{j_s} \epsilon_s 0 \} \right)= \chi(M_F),$$
if $n$ is even and
$$\chi \left( \mathcal{L}_{I} \cap \{f_{j_1} \epsilon_1 0, \ldots, f_{j_s} \epsilon_s 0 \} \right)=2- \chi(M_F),$$
if $n$ is odd, where  for $i \in \{1,\ldots, s \}$, $\epsilon_i \in \{\le,\ge \}$.

\vskip0,2cm

Let us prove first IH(1). We have the following equality:
$$\displaylines{
\qquad \chi \left( \mathcal{L}_{I} \cap \{f_{j_1} \geq 0 \} \right)=
\chi \left(  f_I ^{-1}(0) \cap \{f_{j_1} \ge \delta \}\cap S_\epsilon \right) + \hfill \cr
\hfill \chi \left(  f_I ^{-1}(0) \cap \{ 0 \le f_{j_1} \le \delta \}\cap S_\epsilon \right) -
\chi \left(  f_I ^{-1}(0) \cap \{f_{j_1} = \delta \}\cap S_\epsilon \right) , \qquad \cr
}$$
where $0< \delta \ll \epsilon \ll 1$. As already explained above, by a deformation argument due to Milnor, $f_I^{-1}(0) \cap f_{j_1}^{-1}(\delta) \cap B_\epsilon$ is homeomorphic to $f_I ^{-1}(0) \cap \{f_{j_1} \ge \delta \}\cap S_\epsilon$. So the first term of the right-hand side is equal to $\chi(M_{F})$. By Durfee's result, $ f_I ^{-1}(0) \cap \{ 0 \le f_{j_1} \le \delta \}\cap S_\epsilon$ retracts by deformation to
$ f_I ^{-1}(0) \cap \{  f_{j_1} =0 \}\cap S_\epsilon$ and so the second term is equal to $\chi(\mathcal{L}_{I \cup \{ j_{1}\}}).$
Similarly we have $$\chi(\mathcal{L}_{I}\cap \{f_{j_1}\leq 0\})=\chi(M_{F})+\chi(\mathcal{L}_{I \cup \{ j_{1}\}})-\chi \left(  f_I ^{-1}(0) \cap \{f_{j_1} = -\delta \}\cap S_\epsilon \right).$$
Therefore we can conclude by the previous lemma that
$$\chi \left( \mathcal{L}_{I} \cap \{f_{j_1} \geq 0 \} \right)=\chi \left( \mathcal{L}_{I} \cap \{f_{j_1} \leq 0 \} \right).$$
Applying the Mayer-Vietoris sequence, we can write:
$$
 \chi \left( \mathcal{L}_{I} \right) =\chi \left( \mathcal{L}_{I} \cap \{f_{j_1} \geq 0 \} \right)+
 \chi \left( \mathcal{L}_{I} \cap \{f_{j_1} \leq 0 \} \right)
-\chi \left( \mathcal{L}_{I} \cap \{f_{j_1} = 0 \} \right).
$$
It is easy to conclude using Theorem \ref{CharLink}.

\vspace{0.2cm}

Let us assume that IH($1$),$\ldots$, IH($s-1$) are satisfied and let us prove IH($s$). We have the following equality:
$$\chi \left( \mathcal{L}_{I} \cap \{f_{j_1} \epsilon_1 0, \ldots,f_{j_{s-1}} \epsilon_{s-1} 0, f_{j_s} \geq 0 \} \right)= $$
$$\chi \left(  f_I ^{-1}(0) \cap \{f_{j_1} \epsilon_1 0, \ldots, f_{j_{s-1}} \epsilon_{s-1} 0 , f_{j_s} \ge \delta \} \cap S_\epsilon \right) + $$
$$ \chi \left(  f_I ^{-1}(0) \cap \{f_{j_1} \epsilon_1 0, \ldots, f_{j_{s-1}} \epsilon_{s-1} 0 , 0 \le f_{j_s} \le \delta \}\cap S_\epsilon \right) -$$
$$\chi \left(  f_I ^{-1}(0) \cap\{f_{j_1} \epsilon_1 0, \ldots, f_{j_{s-1}} \epsilon_{s-1} 0 , f_{j_s} = \delta \}\cap S_\epsilon \right) , $$
where $0< \delta \ll \epsilon \ll 1$. By an adaptation of Milnor's deformation argument,  we see that
$$f_I^{-1}(0) \cap  \{f_{j_1} \epsilon_1 0, \ldots, f_{j_{s-1}} \epsilon_{s-1} 0 , f_{j_s} = \delta \} \cap B_\epsilon,$$ is homeomorphic to
$$f_I ^{-1}(0) \cap \{f_{j_1} \epsilon_1 0, \ldots, f_{j_{s-1}} \epsilon_{s-1} 0 , f_{j_s} \ge \delta \} \cap S_\epsilon.$$ So the first term of the right-hand side is equal to $\chi(M_{F})$ by Proposition \ref{CharSemi1}. By Durfee's result,
$$ f_I ^{-1}(0) \cap \{f_{j_1} \epsilon_1 0, \ldots, f_{j_{s-1}} \epsilon_{s-1} 0 , 0 \le f_{j_s} \le \delta \}\cap S_\epsilon,$$ retracts by deformation to
$$ f_I ^{-1}(0) \cap  \{f_{j_1} \epsilon_1 0, \ldots, f_{j_{s-1}} \epsilon_{s-1} 0 , f_{j_s} =0 \}\cap S_\epsilon,$$ and so the second term is equal to $\chi(\mathcal{L}_{I\cup \{j_s\}}\cap \{f_{j_1} \epsilon_1 0, \ldots, f_{j_{s-1}} \epsilon_{s-1} 0 \})$. Applying the previous lemma as for IH(1), we can conclude that
$$\chi \left( \mathcal{L}_{I} \cap \{f_{j_1} \epsilon_1 0, \ldots,f_{j_{s-1}} \epsilon_{s-1} 0, f_{j_s} \geq 0 \} \right)=$$ $$\chi \left( \mathcal{L}_{I} \cap \{f_{j_1} \epsilon_1 0, \ldots,f_{j_{s-1}} \epsilon_{s-1} 0, f_{j_s} \leq 0 \} \right).$$
Applying the Mayer-Vietoris sequence, we can write:
$$ \chi \left( \mathcal{L}_{I} \cap \{f_{j_1} \epsilon_1 0, \ldots, f_{j_{s-1}} \epsilon_{s-1} 0  \}  \right)=$$
$$\chi \left( \mathcal{L}_{I} \cap \{f_{j_1} \epsilon_1 0, \ldots, f_{j_{s-1}} \epsilon_{s-1} 0 , f_{j_s} \ge 0 \}  \right)+ $$
$$ \chi \left( \mathcal{L}_{I} \cap \{f_{j_1} \epsilon_1 0, \ldots, f_{j_{s-1}} \epsilon_{s-1} 0 , f_{j_s} \le 0 \} \right)-$$
$$\chi \left( \mathcal{L}_{I}\cap \{f_{j_1} \epsilon_1 0, \ldots, f_{j_{s-1}} \epsilon_{s-1} 0 , f_{j_s} = 0 \} \right). $$
It is easy to conclude using the induction hypothesis IH($s-1$).   \endproof

\begin{corollary}
If $s \ge 2$ then we have:
$$\chi \left( f_I^{-1}(0)  \cap \{f_{j_1} \epsilon_1 0, \ldots, f_{j_{s-1}} \epsilon_{s-1} 0, f_{j_s} = \delta \} \cap S_\epsilon\right) = \chi(M_F).$$
\end{corollary}
\proof  We use the following equality already mentioned above:
$$\chi \left( \mathcal{L}_{I} \cap \{f_{j_1} \epsilon_1 0, \ldots, f_{j_{s-1}} \epsilon_{s-1} 0 , f_{j_s} \ge 0 \}  \right)= $$
$$\chi \left(  f_I ^{-1}(0) \cap \{f_{j_1} \epsilon_1 0, \ldots, f_{j_{s-1}} \epsilon_{s-1} 0 , f_{j_s} \ge \delta \} \cap S_\epsilon \right) + $$
$$ \chi \left(  f_I ^{-1}(0) \cap \{f_{j_1} \epsilon_1 0, \ldots, f_{j_{s-1}} \epsilon_{s-1} 0 , 0 \le f_{j_s} \le \delta \}\cap S_\epsilon \right) -$$
$$\chi \left(  f_I ^{-1}(0) \cap\{f_{j_1} \epsilon_1 0, \ldots, f_{j_{s-1}} \epsilon_{s-1} 0 , f_{j_s} = \delta \}\cap S_\epsilon \right) , $$
where $0< \delta \ll \epsilon \ll 1$. The term of the left-hand side and the second term of the right-hand side are equal by the previous theorem and the first term of the left-hand side is equal to $\chi(M_F)$, as already explained in the proof of the previous theorem. \endproof

\section{Applications}

In \cite{Sz2} Z. Szafraniec proved interesting formulae relating the Euler number of the link of a weighted homogeneous real polynomial function $f:\mathbb{R}^{n} \to \mathbb{R}$, such that $df(0)=0,$ with the topological degrees of mappings which are explicitly constructed in terms of $f.$ Let us remind the main steps and results.

\vspace{0.2cm}

Let $f:\mathbb{R}^{n} \to \mathbb{R}$ be a weighted homogeneous real polynomial function of type $(d_{1},\cdots,d_{n};d),$ with $df(0)=0,$ and denote by $L=\{x\in S^{n-1}; f(x)=0\}$ the link of $\{f=0\}.$

Let $p$ be the smallest positive integer such that $2p>d$ and each $d_{i}$ divides $p.$ Also denote by $a_{i}=\displaystyle{\frac{p}{d_{i}}}$ and

\begin{equation}
\omega=\frac{x_{1}^{2a_{1}}}{2a_{1}}+\cdots+\frac{x_{n}^{2a_{n}}}{2a_{n}}.
\end{equation}

Now consider $g_{1}=f-\omega$ and $g_{2}=-f-\omega.$

\begin{lemma}[\cite{Sz2}, page 242]

Let $I_{i}=(\frac{\partial g_{i}}{\partial x_{1}},\cdots,\frac{\partial g_{i}}{\partial x_{n}})$ be the ideal generated by the partial derivatives of $g_{i}$, for each $i=1,2,$ in $\mathbb{R}[[x_{1},\cdots,x_{n}]].$ Then, $$\dim_{\mathbb{R}}\frac{\mathbb{R}[[x_{1},\cdots,x_{n}]]}{I_{i}}<\infty .$$

\end{lemma}

In the case where $f$ is homogeneous, i.e., the weight $d_{i}=1$ for all $i=1,\cdots, n$, the author observed that the integer $p=[\frac{d}{2}]+1$ and $$\displaystyle{\omega=(x_{1}^{2p}+\cdots + x_{n}^{2p})/2p}.$$

Now consider the mappings

\begin{equation}
H_{i}=\nabla g_{i}:(\mathbb{R}^{n},0)\to (\mathbb{R}^{n},0),
\end{equation}
for $i=1,2.$ By the lemma above, we have that $0\in \mathbb{R}^{n}$ is isolated in $H_{i}^{-1}(0),$ and so $\deg_{0} H_{i}$ is well defined, for each $i=1,2.$

\vspace{0.2cm}

The next result relates the topological degrees of the mappings $H_{i}$ with the Euler number of the link $L=f^{-1}(0)\cap S^{n-1}.$

\begin{theorem} [\cite{Sz2}, Theorem 5, page 244]

\begin{equation}
\chi(L)=2-(\deg_{0} H_{1} +\deg_{0} H_{2}+\chi(S^{n-1})).
\end{equation}

\end{theorem}

\vspace{0.2cm}

If the degree of homogeneity $d$ is odd, then is possible to do an involution on the sphere $S^{n-1}$ and in this case we have that $\deg_{0} H_{1} =\deg_{0} H_{2} $ and the result above becomes:

\begin{corollary} [\cite{Sz2},Corollary 6, page 244]

\vspace{0.2cm}

If $d$ is odd then

\begin{equation}
\chi(L)=2(1-\deg_{0} H_{1} )-\chi(S^{n-1}).
\end{equation}
\end{corollary}

\vspace{0.2cm}

In what follow we will use Szafraniec's formulae and our previous formulae to compute some examples.

\vspace{0.2cm}

\begin{example}
Consider $f: (\mathbb{C}^{3},0)\to (\mathbb{C},0),$ $f(x,y,z)=x^{2}z+y^{2}.$ It follows from example \ref{R1} that Milnor's conditions $(a)$ and $(b)$ are clearly satisfied. In this case we can apply Sakamoto's formula \cite{Sa} to get that the Milnor fiber $M_{f}$ have the homotopy type of the $2-$dimensional sphere $S^{2},$ and so $\chi(M_{f})=2.$

\vspace{0.2cm}

Let $g=\Re(f):(\mathbb{R}^{6},0)\to (\mathbb{R},0)$ be the  function given by the real part of $f.$ Observe that, since the $\dim (\Sigma_g)>0,$ then the link is not a manifold, therefore we can not find easily the Euler number of the link. Applying our Theorem \ref{CharLink}, we have that $\chi(\mathcal{L}_g)=2\chi(M_{f})=4,$ where $\mathcal{L}_g:=g^{-1}(0)\cap S_{\epsilon}^{5}$ is the link of real function $g.$

\vspace{0.2cm}

Denote $x=x_{1}+ix_{2},~ y=y_{1}+iy_{2}$ and $z=z_{1}+iz_{2}.$  Therefore, the real part $g=z_{1}(x_{1}^{2}-x_{2}^{2})-2z_{2}x_{1}x_{2}+y_{1}^{2}-y_{2}^{2}$ is a weighted-homogeneous polynomial function of type $(2,2,3,3,2,2;6),$ and so we can apply Szafraniec's formula described above as follows:

\vspace{0.2cm}

It is easy to see that $p=6,$ $a_{1}=a_{2}=a_{5}=a_{6}=3,$ $a_{3}=a_{4}=2$ and
$$\omega =\frac{x_{1}^{6}}{6}+\frac{x_{2}^{6}}{6}+\frac{y_{1}^{4}}{4}+\frac{y_{2}^{4}}{4}+\frac{z_{1}^{6}}{6}+\frac{z_{2}^{6}}{6}.$$

\vspace{0.2cm}

It follows from $(3)$ above that the mappings

\vspace{0.2cm}

\begin{center}

$H_{1}=(2x_{1}z_{1}-2z_{2}x_{2}-x_{1}^{5}, -2x_{2}z_{1}-2x_{1}z_{2}-x_{2}^{5}, 2y_{1}-y_{1}^{3},-2y_{2}-y_{2}^{3},x_{1}^{2}-x_{2}^{2}-z_{1}^{5},-2x_{1}x_{2}-z_{2}^{5}),$
\end{center}

\vspace{0.2cm}

\begin{center}
$H_{2}=(-2x_{1}z_{1}+2z_{2}x_{2}-x_{1}^{5}, 2x_{2}z_{1}+2x_{1}z_{2}-x_{2}^{5}, -2y_{1}-y_{1}^{3},2y_{2}-y_{2}^{3},-x_{1}^{2}+x_{2}^{2}-z_{1}^{5},2x_{1}x_{2}-z_{2}^{5}).$
\end{center}

Now, computations shows that $\deg_{0} H_{j} =-1,~ j=1,2$ and by Szafraniec's formula $(6)$ we have that

$$\chi(\mathcal{L}_{g})=2-(-1-1+0)=4,$$
and so, both results coincide.

\end{example}

\begin{example}
Let $f=(P,Q): (\mathbb{R}^{3},0)\to (\mathbb{R}^{2},0)$ $f(x,y,z)=(zx^2+zy^2+y^3,x).$ It is easy to see that $\Sigma_f=\{(0,0,z); z\in \mathbb{R}\},$ $\Sigma_f\subset V,$ and so Milnor's condition $(a)$ holds.


On the one hand we have that $Q=x$ and since its link is diffeomorphic to $S^{1}$ then $\chi(\mathcal{L_{Q}})=0.$

\vspace{0.2cm}

On the other hand, since $P=zx^2+zy^2+y^3$ is homogeneous of degree $d=3,$ then $p=2,$ $g_{1}=zx^2+zy^2+y^3-\omega ,$ where $\omega=\frac{x^{4}}{4}+\frac{y^{4}}{4}+\frac{z^{4}}{4}$ and
$$H_{1}=(2xz-x^3,3y^2+2yz-y^3,x^2+y^2-z^3).$$

It is easy to see that $\deg_{0} H_{1} =-1$ and by formula $(7)$ above we have $\chi(\mathcal{L}_{P})=2.$ Therefore, Milnor's condition $(b)$ cannot be satisfied.

\vspace{0.3cm}


\end{example}

\begin{example} The next example comes from \cite{TYA}, Example 5.1, Section 5. Applying a Thom-Sebastiani type Theorem as explained in Section 5 of cited paper, it is possible to produce examples in all odd dimension.

\vspace{0.3cm}

Let $f=(P,Q): (\mathbb{R}^{5},0)\to (\mathbb{R}^{2},0),$ $f(x,y,z,u,v)=(y^{4}-z^2x^2-x^4+u^2-v^2,xy+2uv).$ It is easy to see that $\Sigma_{f}\subseteq V.$ It was proved in \cite{TYA} that Milnor's condition $(b)$ follows as an application of the Curve Selection Lemma.

\vspace{0.2cm}

Calculations shows that $\deg_{0} {H_{i}} =1,$ for $i=1,2.$ So, $\chi(L_{Q})=-2,$ $\chi(M_{f})=2$ and $\chi(\partial M_{f})=4.$ Therefore, we can claim that the boundary of the Milnor fiber is not connected, since it contains at least two disjoint copies of $S^2.$

\end{example}

\begin{remark}

Let $f_{1},\cdots , f_{s}:U\subseteq \mathbb{R}^{n}\to \mathbb{R}$ be analytic functions with $f_{1}(0)=\cdots f_{s}(0)=0$ and let $f(x)=f_{1}(x)^2+\cdots +f_{s}(x)^2.$ The following equality of analytic sets holds $$\{x\in U;~ f_{1}(x)=\cdots =f_{s}(x)=0\}=\{x\in U;~ f(x)=0\}.$$

In \cite{Sz1} the author considered the function $g(x)=f(x)-c(x_{1}^{2}+\cdots +x_{n}^{2})^{k},$ where $c>0$ and $k$ an integer, and showed that for $k$ large enough the function $g$ has an isolated singular point at the origin. Moreover, he proved that for all small radius $\epsilon $, the following Poincar\'e-Hopf type formula holds true: $$\chi(\mathcal{L}_{f})=\chi(\{x\in S_{\epsilon}^{n-1}; ~ f(x)=0\})=1-\deg_{0} \nabla g.$$

\vspace{0.2cm}

Therefore, connecting this result with our previous formulae, we can conclude that for a given analytic mapping satisfying Milnor's conditions $(a)$ and $(b)$ the Euler-Poincar\'e number of the Milnor fiber also satisfies a Poincar\'e-Hopf type formula.

\end{remark}

\end{document}